\documentclass[12pt,twoside,doublespace]{article}

\usepackage{amssymb}
\usepackage{graphicx}

\def\smskip{\par\vskip 5 pt}
\def\QED{\hfill $\Box$\smskip}
\newtheorem{theorem}{Theorem}

\newtheorem{proposition}{Proposition}

\newtheorem{example}{Example}
\newtheorem{remark}{Remark}

\begin{document}

\begin{center}

\vspace{35pt}

{\Large \bf A General Class of }

\vspace{5pt}

{\Large \bf Relative Optimization Problems

}

\vspace{35pt}

{\sc I.V.~Konnov\footnote{\normalsize Department of System Analysis
and Information Technologies, Kazan Federal University, ul.
Kremlevskaya, 18, Kazan 420008, Russia.\\ E-mail: konn-igor@ya.ru}}

\end{center}

\begin{abstract}
We consider relative or subjective optimization problems where
the goal function and feasible set are dependent of the current state
of the system under consideration. In general,
 they are formulated  as quasi-equilibrium problems, hence
 finding their solutions may be rather difficult.
 We describe a rather general class of relative
optimization problems in metric spaces, which in addition depend on the starting state.
 We also utilize quasi-equilibrium type formulations of
these problems and show that they admit rather simple descent solution methods.
This approach gives suitable trajectories tending to a relatively
optimal state. We describe several examples of applications of these problems.

\end{abstract}


\section{Introduction}\label{sc:1}

The usual requirement to choose the best variant in various decision making problems
naturally leads to their optimization formulations. That is, one then has to find an element
attributed to a decision from some given feasible set $D$ that yields
the maximal (or minimal) value of some goal (utility) function $\varphi$.
For brevity, we write this problem as
\begin{equation} \label{eq:1.1}
 \max \limits _{y \in D} \to \varphi(y).
\end{equation}
However, due to incomplete and inexact knowledge about
the goal function and feasible set this simple formulation
usually needs certain corrections; see e.g. \cite{HCB04,BEN09}.
Recently, a new approach to this problem was proposed in \cite{Kon19d}
where it was supposed that the presentation of the goal and constraints defining the system model
may vary together with the changes of the system state and that only some limited information
about the goal and  constraints  may be known at each state.
It was proposed to consider such mathematical models as relative or subjective
optimization problems with respect to system states and
 to formulate them as (quasi-)equilibrium problems.
 This means that the goal function is replaced with a bi-function $\phi(x,y)$
 so that $\phi(x,\cdot)$ is the goal function attributed to a current state $x$.
Similarly, the feasible set $D$ may also depend on the
states and is replaced with a set-valued mapping $x \mapsto D (x)$.
That is, we have only restricted knowledge about the problem at each point.
A relatively optimal state $x^{*}$ should give the maximal value of the goal
function which is compared with all the other feasible states evaluated at
the current state $x^{*}$, i.e. one has to find
$x^{*} \in D(x^{*})$ such that
\begin{equation} \label{eq:1.1}
\phi(x^{*},x^{*}) \geq \phi(x^{*},y) \quad \forall y \in D(x^{*}).
\end{equation}
It follows that the above concept gives certain restricted optimality.
Nevertheless, it can be used in order to decide whether the current state
is suitable or should be changed, thus implementing a weaker solution concept.
We observe that (\ref{eq:1.1}) is nothing but the so-called {\em quasi-equilibrium problem}
(QEP for short); see \cite{BL84,Har91,Aub98}.
Finding a solution of quasi-equilibrium problems may be rather difficult because of the
presence of the moving feasible set.

In this paper, we describe a rather general class of relative
optimization problems, which in addition depend on the starting state.
 We also take  quasi-equilibrium type formulations of
these problems and propose simple descent solution methods for creating
suitable trajectories to a relatively optimal state. We establish
existence results for these problems under mild conditions and give
illustrative examples of applications.


\section{Basic Problem Formulations}\label{sc:2}

We first describe a general model of a system whose possible states
are contained in a set $X \subseteq E$ where $E$ is a metric space.
The starting state $x^{0} \in X$ is known.
Given a state $x \in X$, one can define the set of
feasible states $D(x)$. This means that the system can move from $x$ to any
$y \in D(x)$ and the utility estimate $\varphi(x,y)$ of any state
$y \in D(x)$ is known at $x$, i.e. $D(x)$ stands for a \lq\lq trust region" at $x$.
We suppose that the estimate $u(x)= \varphi(x,x)$
is precise, but the value $u(y)= \varphi(y,y)$ is not supposed to be known at $x$.
It follows that $x \in D(x)$ for any $x \in X$. Next, each move $(x \to y)$
requires certain expenses $c(x,y)$. We suppose that $c(x,y)$ is non-negative and
known at $x$ for any $y \in D(x)$. Hence, we can define the estimate of
pure expenses for the move $(x \to y)$ as follows
$$
f(x,y)=\varphi(x,x)+c(x,y)-\varphi(x,y),
$$
as well as the precise pure expenses for this move
$$
e(x,y)=u(x)+c(x,y)-u(y).
$$
Choice of the set $D(x)$ at $x \in X$ should guarantee that the estimates
have some sufficient precision.
We will say that a sequence $\{x^{k}\} \subset X$ is a {\em feasible trajectory}
if $x^{k+1} \in D(x^{k})$ for each number $k$.
Then we can define two relative optimization problems.

\medskip
\noindent {\bf Problem (P1)} {\em Find a point $x^{*} \in X$ such that
\begin{equation} \label{eq:2.1}
f(x^{*},y) \geq 0 \quad \forall y \in D(x^{*}).
\end{equation}}

\medskip

\noindent {\bf Problem (P2)} {\em Find a feasible trajectory $\{x^{k}\}$
with the initial state $x^{0} \in X$ and non-positive pure expenses estimates
such that it either terminates at a solution of Problem  {\bf (P1)} or  its limit points
are solutions of Problem  {\bf (P1)}.}
\medskip

It is clear that (\ref{eq:2.1}) coincides with (\ref{eq:1.1}) if $c(x,x)=0$ and we set
$$
\phi(x,y)=\varphi(x,y)-c(x,y).
$$
We observe that Problem  {\bf (P1)} is stationary since it does not depend on the initial state
whereas Problem  {\bf (P2)} depends on the initial state essentially.
In fact, then one also has to take a feasible trajectory $\{x^{k}\}$ such that
$f(x^{k-1},x^{k}) \leq 0$ for each $k$. Then we have
$$
f(x^{0},x^{1})+f(x^{1},x^{2})+ \ldots +f(x^{k-1},x^{k}) \leq 0
$$
for each $k$, i.e. we intend to move the system from the current state to a
relatively optimal state without expenses. Existence of a  solution of Problem  {\bf (P2)}
means that sequential taking some moves without expenses can yield a
relatively optimal state. It should be noted that Problem  {\bf (P2)} differs from
the usual global discrete time  optimal control problems; see e.g. \cite{Zas06}.

\begin{remark} \label{rm:2.1}
We note that the usual decision making approaches require the choice of the best variant
with respect to some given optimality criterion even in the presence of uncertainty factors.
That is, such a solution must be globally optimal with respect to all the variants.
However, we think that the \lq\lq globally marginal" behaviour is not so suitable in the case of
inexact and incomplete data. The above relaxed optimality concepts give an alternative
approach, which enables one to only evaluate the necessity to change the current state
of the system. This means that the optimization formulation is then restricted within
a variable feasible set containing only the states whose estimates at the current state
are sufficiently precise.
\end{remark}


\section{The Basic Method and Its Convergence} \label{sc:3}

We will use the following set of basic assumptions.

\noindent
{\bf (A1)}  {\em The set $X \subseteq E$ is nonempty and closed, the bi-function
$\varphi : X \times X \rightarrow \mathbb{R}$ is continuous, the bi-function
$c : X \times X  \rightarrow \mathbb{R}$ is non-negative and continuous,
and $c(x,x)=0$ for each $x \in X$.}

\noindent
{\bf (A2)} {\em For  any number $\alpha $ the set
$$
X_{\alpha}=\left\{ x \in X \ | \ u (x) \geq \alpha \right\}
$$
is compact, for  any bounded set $\tilde X \subset X$ there exists a number $\beta $ such that
$\tilde X \subseteq X_{\beta}$.}

\noindent
{\bf (A3)} {\em The mapping $D : X \rightarrow \Pi (X)$
is lower semi-continuous on  $X$ and $x \in D(x)$ for each $x \in X$.}

We recall that a set-valued mapping $T : E \rightarrow \Pi (E)$ is said to be
{\em  lower semi-continuous} at a point  $z \in X$ on a set $X$ if, for any
sequence $\{ x^{k}\} \to z$, $ x^{k} \in X$, and any $t \in T(z)$ there exists
a sequence $\{ t^{k}\} \to t$, $ t^{k} \in T(x^{k})$.
The mapping $T$ is said to be  lower semi-continuous on the set $X$
if it is  lower semi-continuous at any point of $ X$. Here $\Pi(A)$ denotes the family
of all subsets of a set $A$.

Clearly, {\bf (A2)} is a general coercivity condition, which implies that
the usual optimization problem
\begin{equation} \label{eq:3.1}
 \max \limits _{x \in X} \to u(x)
\end{equation}
has a solution and that
$$
u^{*}= \max \limits _{x \in X}  u(x) < +\infty.
$$
We now describe a general threshold descent method (TDM) for Problem  {\bf (P2)} and hence for {\bf (P1)} as well.

\medskip
\noindent {\bf Method (TDM).}
Take the given point $x^{0}$, choose a sequence
$\{\delta _{l}\} \searrow 0$. Set  $l=1$, $k=0$, $ z^{0} =x^{0}$.

For each $k=0,1,\ldots$, we have a point $z^{k}\in X$.
Find a point $z^{k+1} \in D(z^{k})$ such that
\begin{equation} \label{eq:3.2}
f(z^{k},z^{k+1}) < -\delta _{l}.
\end{equation}
If this point does not exist, set $ x^{l} =z^{k}$,
$l= l+1$. Go to the $(k+1)$-th iteration.
\medskip

Therefore, $\delta _{l}$ stands for the current descent threshold, which determines the sufficient profit
for the movement.

In order to guarantee convergence of (TDM)
we need additional conditions for the accuracy of utility estimates
related to system moves expenses. For brevity, set $[\alpha]_{+}=\max \{\alpha, 0\}$ for a number $\alpha$ and
$$
b(x,y)=[\varphi(x,y)-u(y)]_{+}.
$$
That is, $b(x,y)$ is the utility over-estimate of the state $y$ at $x$.

\medskip
\noindent
{\bf (A4)} \hfill
\begin{enumerate}
   \item[{(i)}]    {\em  For  any feasible trajectory $\{z^{k}\}$ it holds that
$$
\lim \limits_{k \to \infty } [b(z^{k},z^{k+1})-c(z^{k},z^{k+1})]_{+} =0;
$$}
   \item[{(ii)}]  {\em  For  any unbounded feasible trajectory $\{z^{k}\}$ it holds that
$$
\sum \limits_{k=0}^{\infty } [b(z^{k},z^{k+1})-c(z^{k},z^{k+1})]_{+} <\infty.
$$}
\end{enumerate}


\begin{theorem} \label{thm:3.1}
Let assumptions {\bf (A1)}--{\bf (A4)} be fulfilled. Then
the sequence $\{x^l \}$ generated by Method (TDM) has limit points,  all these limit points
are solutions of Problem  {\bf (P1)}, and the sequence $\{z^{k}\}$ solves Problem  {\bf (P2)}.
\end{theorem}
{\bf Proof.}
The assertion will be proved in several steps.

\textit{Step 1: For each $l$ the number of changes of the index $k$ is finite.} \\
From the definitions and (\ref{eq:3.2}) we have
\begin{eqnarray*}
 \displaystyle e(z^{k},z^{k+1}) &=& f(z^{k},z^{k+1}) + (\varphi(z^{k},z^{k+1}) - u(z^{k+1})) \leq f(z^{k},z^{k+1}) + b(z^{k},z^{k+1})  \\
 \displaystyle
 &<& -\delta _{l}+ b(z^{k},z^{k+1}),
\end{eqnarray*}
hence
\begin{equation} \label{eq:3.3}
u(z^{k+1})-u(z^{k}) > \delta _{l}-[b(z^{k},z^{k+1})-c(z^{k},z^{k+1})]_{+}
\end{equation}
for each fixed index $l$.
If the number of changes of the index $k$ is infinite for some $l$, {\bf (A4)} (i) and (\ref{eq:3.3}) imply
$u(z^{k}) \to +\infty$ as $k \to \infty$, which is a contradiction.

\textit{Step 2:  The sequence $\{z^{k}\}$ is bounded.} \\ Suppose
$\{z^{k}\}$ is unbounded. Then {\bf (A4)} (ii) and (\ref{eq:3.3}) imply that
$$
\lim \limits_{k \to \infty } u(z^{k})=\tilde u < +\infty
$$
due to Lemma 1 in \cite[Chapter III]{GT89}. Hence, there exist numbers $\alpha$ and $k'$ such that
$z^{k} \in X_{\alpha}$ if $k > k'$. It follows that the sequence $\{z^{k}\}$ is contained in
the compact set $X_{\alpha}$, which is a contradiction.

\textit{Step 3:  The sequence $\{x^l \}$ has limit points,  all these limit points
are solutions of Problem  {\bf (P1)}.} \\ From Steps 1--2 it follows that
the sequence $\{x^l \}$ is infinite and bounded, hence it
is contained in a compact set $X_{\beta}$ due to {\bf (A2)}. It follows that
$\{x^l \}$ has limit points.
For each $l$ from the definition we have
\begin{equation} \label{eq:3.4}
f(x^l,y) \geq -\delta _{l} \quad \forall y \in D(x^l).
\end{equation}
Let $\bar x$ be an arbitrary limit point of
$\{x^{l}\}$, i.e. $ \{x^{l_{s}} \} \to \bar x$. Then $\bar x \in X$ since
$X$ is closed. Take any
$\bar y \in D(\bar x)$, then there exists a sequence of points  $ \{ y^{l_{s}} \}$, $\{
y^{l_{s}} \} \to \bar y$ such that $ y^{l_{s}} \in D(x^{l_{s}})$ since
the mapping $D$ is lower semi-continuous on  $X$.
Setting $l=l_{s}$ and $y=y^{l_{s}}$ in (\ref{eq:3.4}) and
taking the limit $s \to \infty$ give
$$
f(\bar x,\bar y) \geq 0,
$$
i.e. $\bar x$ is a solution of Problem  {\bf (P1)}. Since $f(z^{k},z^{k+1})<0$,
$\{z^{k}\}$ is a solution of Problem  {\bf (P2)}.
\QED

Clearly, Theorem \ref{thm:3.1} implies existence of solutions of Problems  {\bf (P1)} and {\bf (P2)} under
assumptions {\bf (A1)}--{\bf (A4)}. We observe that a solution of the optimization problem (\ref{eq:3.1})
is not in general a solution of Problem  {\bf (P1)} under
assumptions {\bf (A1)}--{\bf (A4)}  as the following simple examples illustrate.


\begin{example} \label{ex:3.1}
{\em Let $X=[0,1]$, $u(x)=1-x/4$, $c(x,y)\equiv 0$, $D(x)=[x,x+0.1(1-x)]\cap X$,
$\varphi(x,y)=(1-y/4)+0.6|0.5-x|(y-x)$. Hence  $\varphi(x,y)=u(y)$ if
$x= 0.5$. Then the point  $x^{0}=0$ is a unique solution of (\ref{eq:3.1}) since $u(x^{0})=1$.
But it is not a solution of {\bf (P1)} since $x^{1}=0.1 \in D(x^{0})$ and
$$
f(x^{0},x^{1}) =x^{1}(0.25-0.3)< 0.
$$
The point $\bar x=1/12$ is the solution of {\bf (P1)} closest to $x^{0}$ since
$$
f(\bar x,y) =0 \quad \forall y \in D(\bar x).
$$
At the same time, we conclude that all the assumptions in {\bf (A1)}--{\bf (A4)} are fulfilled.
In fact, any feasible trajectory $\{z^{k}\} $ is bounded and $z^{k} \leq z^{k+1}$.
Hence it converges to a point in  $X$, which implies
$$
\lim \limits_{k \to \infty } b(z^{k},z^{k+1}) =0.
$$}
\end{example}


\begin{example} \label{ex:3.2}
{\em Let $X=[0,1]$, $u(x)=1-x/4$, $c(x,y)\equiv 0$,
\begin{eqnarray*}
 \displaystyle \varphi(x,y) &=& (1-y/4)+[0.5-x]_{+}(y-x),  \\
 \displaystyle D(x) &=& [x-0.1[0.5-x]_{+},x+0.1(2-x)]\cap X.
\end{eqnarray*}
Here $D(x)$ is not a singleton at any point $x \in X$ and
  $\varphi(x,y)=u(y)$ if $x \geq 0.5$. Again the point  $x^{0}=0$ is a unique solution
of (\ref{eq:3.1}) since $u(x^{0})=1$.
But it is not a solution of {\bf (P1)} since $x^{1}=0.2 \in D(x^{0})$ and
$$
f(x^{0},x^{1}) =x^{1}(0.25-0.5)< 0.
$$
The point $\bar x=0.25$ is the solution of {\bf (P1)} closest to $x^{0}$ since
$$
f(\bar x,y) =0 \quad \forall y \in D(\bar x).
$$
Also, all the assumptions in {\bf (A1)}--{\bf (A4)} are fulfilled.
It suffices to check {\bf (A4)} (i). Let us take
any feasible trajectory $\{z^{k}\} $. If $z^{k} \leq 0.5$, then $z^{k} \leq z^{k+1}$, but if
$z^{k} \geq 0.5$, then $z^{k+1} \geq 0.5$ and  $b(z^{k},z^{k+1}) =0$.
It follows that only one transition $(z^{k} \leq 0.5) \to (z^{k+1} > 0.5) $
is possible for $\{z^{k}\} $ and that   {\bf (A4)} (i) is fulfilled.}
\end{example}


\section{Discussion of Conditions and Modifications} \label{sc:4}

We observe that conditions {\bf (A1)}--{\bf (A3)} seem rather natural and simple.
They even do not involve convexity/ monotonicity properties and
do not impose restrictions on the values of the mapping $x \mapsto D (x)$.
Therefore, the set of assumptions is somewhat different from the custom ones; cf. e.g.
\cite{BL84,YT97,Aub98}.
We now discuss the assumptions in {\bf (A4)} which in fact indicate the precision bounds for
utility estimates of any state $y \in D(x)$ at $x$. In the general case
the cost value $c (z^{k},z^{k+1})$ is known at $z^{k}$ by assumption. Hence,
the proper choice of the set $D(z^{k})$ needs certain concordance of the utility over-estimate
and move expenses for providing the relation
$$
[b(z^{k},z^{k+1})-c(z^{k},z^{k+1})]_{+} \approx 0
$$
and attaining the convergence. In other words, the difference between the utility over-estimate
and move expenses should tend to zero along any infinite feasible trajectory
and this convergence should be rather rapid if the trajectory is unbounded.

Let us take the modified pair of conditions.

\medskip
\noindent
{\bf (A2${}'$)} {\em For  any number $\alpha $ the set
$$
X_{\alpha}=\left\{ x \in X \ | \ u (x) \geq \alpha \right\}
$$
is compact.}

\medskip
\noindent
{\bf (A4${}'$)}  {\em  For  any feasible trajectory $\{z^{k}\} $  it holds that
$$
\sum \limits_{k=0}^{\infty } [b(z^{k},z^{k+1})-c(z^{k},z^{k+1})]_{+} <\infty.
$$}
\medskip

The assertions of Theorem \ref{thm:3.1} remain true if we replace {\bf (A2)} and {\bf (A4)} with
{\bf (A2${}'$)} and {\bf (A4${}'$)}, respectively. Here {\bf (A2${}'$)} is weaker than {\bf (A2)},
but  {\bf (A4${}'$)} is stronger than {\bf (A4)}. Nevertheless, this is the case if the utility over-estimate of a state $y\in D(x)$ at $x$
appears to be less than the move expenses $c(x,y)$ due to our subjective choice of the set $D(x)$.
Then we can in turn replace {\bf (A2${}'$)}  and {\bf (A4${}'$)} with the following.

\medskip
\noindent
{\bf (A2${}''$)}  {\em For  some number $\alpha \leq u (x^{0})$ the set
$X_{\alpha}$ is compact.}
\medskip

\medskip
\noindent
{\bf (A4${}''$)}  {\em  For any $x \in X$  it holds that
$$
b(x,y)\leq c(x,y) \quad \forall y \in D(x).
$$}
\medskip

The assertions of Theorem \ref{thm:3.1} remain true if we replace {\bf (A2)} and {\bf (A4)} with
{\bf (A2${}''$)} and {\bf (A4${}''$)}, respectively.
Let us now suppose that the cost bi-function
$c $ satisfies {\bf (A1)} without any additional assumptions. Then {\bf (A4)} should be modified as follows.

\medskip
\noindent
{\bf (A5)} \hfill
\begin{enumerate}
   \item[{(i)}]    {\em  For  any feasible trajectory $\{z^{k}\} $  it holds that
$$
\lim \limits_{k \to \infty } b(z^{k},z^{k+1}) =0;
$$}
   \item[{(ii)}]  {\em  For  any unbounded feasible trajectory $\{z^{k}\} $  it holds that
$$
\sum \limits_{k=0}^{\infty } b(z^{k},z^{k+1}) <\infty.
$$}
\end{enumerate}
\medskip

This means that only the utility over-estimates tend to zero along any infinite feasible  trajectory
and that this convergence is rather rapid if the trajectory is unbounded. This property can be
invoked by the usual training process along the trajectory and by the proper choice of the sets $D(x)$.
As above, we can use proper modifications of {\bf (A5)} by analogy with  {\bf (A4${}'$)} and {\bf (A4${}''$)}.
For instance,  the assumptions in {\bf (A5)} clearly hold true if there is no any over-estimate, i.e. when
$\varphi(x,y) \leq u(y)$ for any $y \in D(x)$. In this case {\bf (A2)} can be replaced with {\bf (A2${}''$)}.
Then the assertions of Theorem \ref{thm:3.1} remain true.

Let us take the simple descent method (SDM) for Problem  {\bf (P2)}:
\begin{equation} \label{eq:4.1}
x^{k+1} \in D(x^{k}), \ f(x^{k},x^{k+1}) < 0 \quad  \mbox{for} \ k= 0,1,\ldots
\end{equation}
Unlike (TDM), it does not converge to a solution under more strong assumptions
as the following simple example illustrates.


\begin{example} \label{ex:4.1}
{\em Let $X=[0,1]$, $u(x)=x$, $c(x,y)\equiv 0$, $D(x)\equiv X$. Then the process
$$
x^{k+1}=x^{k}+2^{-(k+2)}, \ k= 0,1,\ldots, \ x^{0}=0,
$$
which corresponds to (\ref{eq:4.1}), clearly converges to $\tilde x=0.5$
instead of the unique solution $x^{*}=1$. }
\end{example}
However, (SDM)  can be useful in the
case where the set $X$ is countable and there exists a lower positive threshold for
move expenses. Then we can remove all the continuity assumptions and
 modify the conditions in {\bf (A1)}--{\bf (A4)} as follows.

\medskip

\noindent
{\bf (B1)}  {\em The set $X \subseteq E$ is nonempty and countable, $x \in D(x)$ for each $x \in X$,
$c(x,x)=0$ for each $x \in X$,  and there exists a number $\delta>0$ such that
$c(x,y)\geq \delta$ for all $x,y \in X$, $x\neq y$.}

\medskip

\noindent
{\bf (B2)} \hfill
\begin{enumerate}
   \item[{(i)}]  {\em   It holds that
   $$
u^{*}= \sup \limits _{x \in X}  u(x) < +\infty;
$$}
   \item[{(ii)}] {\em For  any feasible trajectory $\{z^{k}\}$  it holds that
$$
\lim \limits_{k \to \infty } b(z^{k},z^{k+1}) =0.
$$}
\end{enumerate}


\begin{proposition} \label{pro:4.1}
Let assumptions {\bf (B1)}--{\bf (B2)} be fulfilled. Then
the sequence $\{z^{k} \}$ generated by Method (SDM) solves Problem  {\bf (P2)}.
It is finite and stops at a solution of Problem  {\bf (P1)}.
\end{proposition}
{\bf Proof.}
It suffices to prove the finiteness of Method (SDM).
From the definitions and (\ref{eq:3.2}) we  have
\begin{eqnarray*}
 \displaystyle e(z^{k},z^{k+1}) &=& f(z^{k},z^{k+1}) + (\varphi(z^{k},z^{k+1}) - u(z^{k+1})) \leq f(z^{k},z^{k+1}) + b(z^{k},z^{k+1})  \\
 \displaystyle
 &<&  b(z^{k},z^{k+1}),
\end{eqnarray*}
hence
$$
u(z^{k+1})-u(z^{k}) > \delta -b(z^{k},z^{k+1})
$$
for each fixed index $l$. If the sequence $\{z^{k}\}$ is infinite, {\bf (B2)} (ii) now implies
$u(x^{k}) \to +\infty$ as $k \to \infty$, which is a contradiction with {\bf (B2)} (i).
\QED

The basic assumptions can be modified in a complete metric space setting.
Then we can remove the compactness assumption.

\medskip

\noindent
{\bf (C1)}  {\em The set $X \subseteq E$ is nonempty and closed,
$E$ is a complete metric space with the metric bi-function
$d : X \times X  \rightarrow \mathbb{R}$.}

\medskip

\noindent
{\bf (C2)}  {\em The bi-functions
$\varphi : X \times X \rightarrow \mathbb{R}$
and $c : X \times X  \rightarrow \mathbb{R}$ are continuous,
$$
u^{*}= \sup \limits _{x \in X}  u(x) < +\infty.
$$}

\medskip

\noindent
{\bf (C3)}  {\em The bi-function
$c : X \times X  \rightarrow \mathbb{R}$ satisfies the triangle inequality, i.e.,
$$
c(x,z)+c(z,y) \geq c(x,y) \quad \forall x,y,z \in X;
$$
there exists an increasing continuous function $\theta
:\mathbb{R} \rightarrow \mathbb{R}$ such that $\theta (0)=0$ and that for all $x,y
\in X$ we have $\theta [d (x,y)] \leq c (x,y)$.}
\medskip

\noindent
{\bf (C4)} {\em For any feasible trajectory $\{z^{k}\}$ it holds that
$$
\sum \limits_{k=0}^{\infty } b(z^{k},z^{k+1}) <\infty.
$$}

\medskip


\begin{theorem} \label{thm:4.1}
Let assumptions {\bf (C1)}--{\bf (C4)} and {\bf (A3)} be fulfilled. Then
the sequence $\{x^l \}$ generated by Method (TDM) converges to
a solution of Problem   {\bf (P1)}, and the sequence $\{z^{k}\}$ solves Problem  {\bf (P2)}.
\end{theorem}
{\bf Proof.}
The assertion will be proved in several steps.

\textit{Step 1: For each $l$ the number of changes of the index $k$ is finite.} \\
From the definitions and (\ref{eq:3.2}) we have
\begin{eqnarray*}
 \displaystyle e(z^{k},z^{k+1}) &=& f(z^{k},z^{k+1}) + (\varphi(z^{k},z^{k+1}) - u(z^{k+1})) \leq f(z^{k},z^{k+1}) + b(z^{k},z^{k+1})  \\
 \displaystyle
 &<& -\delta _{l}+ b(z^{k},z^{k+1}),
\end{eqnarray*}
hence
\begin{equation} \label{eq:4.2}
u(z^{k+1})-u(z^{k}) > \delta _{l}+c(z^{k},z^{k+1})-b(z^{k},z^{k+1})\geq \delta _{l}-b(z^{k},z^{k+1})
\end{equation}
for each fixed index $l$. If the number of changes of the index $k$ is
infinite for some $l$, {\bf (C4)} and (\ref{eq:4.2}) imply
$u(z^{k}) \to +\infty$ as $k \to \infty$, which contradicts {\bf (C2)}.

\textit{Step 2:  The sequence $\{z^{k}\}$ converges to
a point $\bar x \in X$.} \\ From {\bf (C4)} and (\ref{eq:4.2}) we have
\begin{equation} \label{eq:4.3}
\lim \limits_{k \to \infty } u(z^{k})=\tilde u < +\infty
\end{equation}
due to Lemma 1 in \cite[Chapter III]{GT89}. It also follows from (\ref{eq:4.2}) that
$$
c(z^{k},z^{k+1}) \leq u(z^{k+1})-u(z^{k})+b(z^{k},z^{k+1}).
$$
Take any indices $k$ and $m=k+p$, then we have
\begin{eqnarray*}
 \displaystyle \theta [d (z^{k},z^{k+p}) ] &\leq & c(z^{k},z^{k+p}) \leq c(z^{k},z^{k+1})+ \ldots +c(z^{k+p-1},z^{k+p})  \\
 \displaystyle &\leq & u(z^{k+p})-u(z^{k})+\sum \limits_{i=k}^{k+p-1 } b(z^{i},z^{i+1}).
\end{eqnarray*}
On account of {\bf (C3)}, {\bf (C4)} and (\ref{eq:4.3}) we now obtain that for any number $\alpha>0$ there
exists an index $k'$ such that $d (z^{k},z^{m}) < \alpha$ if $\min\{k,m\}> k'$.
Hence, $\{z^{k}\}$ is a Cauchy sequence and it  converges to
a point $\bar x \in X$ since $X$ is closed.

\textit{Step 3:  The sequence $\{x^l \}$ converges to
a point $\bar x \in X$,  which is a solution of Problem  {\bf (P1)}.} \\
Since the sequence $\{x^l \}$ is contained in $\{z^{k}\}$ and is
infinite due to Step 1, Step 2 implies that $\{x^l \}$ converges to
a point $\bar x \in X$.
For each $l$ from the definition we have
\begin{equation} \label{eq:4.4}
f(x^l,y) \geq -\delta _{l} \quad \forall y \in D(x^l).
\end{equation}
Take any $\bar y \in D(\bar x)$, then by {\bf (A3)} there exists a sequence of points  $ \{ y^{l} \}$, $\{
y^{l} \} \to \bar y$ such that $ y^{l} \in D(x^{l})$ since
the mapping $D$ is lower semi-continuous on  $X$.
Setting $y=y^{l}$ in (\ref{eq:4.4}) and
taking the limit $l \to \infty$ give
$f(\bar x,\bar y) \geq 0$,
i.e. $\bar x$ is a solution of Problem  {\bf (P1)}. Since $f(z^{k},z^{k+1})<0$,
$\{z^{k}\}$  is a solution of Problem  {\bf (P2)}.
\QED

\begin{remark} \label{rm:4.1}
The basic technique for obtaining the assertion of Step 2 of Theorem \ref{thm:4.1}
resembles that of the Caristi fixed point theorem; see e.g. \cite[Section 1.8]{Aub98}.
However, $c$ need not be a metric bi-function, besides, we do not determine a choice mapping,
since the  mapping $D$ only imposes restrictions on the choice at a current point,
which should also conform to the descent rule.
For this reason, the set of assumptions is somewhat different.
\end{remark}


\section{Examples of Models} \label{sc:5}

We now describe some applied models, which can be formulated within the proposed
framework. These models are modifications and extensions of those from
\cite{Kon19d,KMT98}.


\begin{example} \label{ex:5.1} {\bf (Treatment of industrial wastes).}
{\em Let us consider an industrial firm which may utilize $n$
production technologies and have a plant for treatment of
its wastes containing $m$ polluted substances. Let
$x=(x_{1},\ldots,x_{n})^{\top} \in \mathbb{R}^{n}$ be the
vector of technology activity levels (activity profile) of the firm. Then $q(x)=
(q_{1}(x),\ldots,q_{m}(x))^{\top} \in \mathbb{R}^{m}$ is
the corresponding vector of its wastes and $\mu(x)$ is the
benefit of this firm. That is, $\mu(x)=\mu_{1}(x)-\mu_{2}(x)$,
where $\mu_{1}(x)$ is the income from  selling its products
and $\mu_{2}(x)$ is the total resource expenses at the technology
activity profile $x$. We denote by $X \subseteq
\mathbb{R}^{n}_{+}$  the whole feasible activity profile set of the firm,
which stands for the set of feasible states.

Next, suppose that the vector $p$ of unit treatment charges
depends on the pollution volumes, that is $p=p[q(x)]$, but
the exact values of these parameters are not known.
Namely, if $x$ is the current
vector of activity levels, then one can calculate the values of the functions
$p_{i}[q(y)]$ only if $q(y)$ belongs to some neighborhood $U(q)$ of $q=q(x)$,
i.e. we have in fact $p_{i}=p_{i}[q(x),q(y)]$. That is, the utility (profit) value estimate at
$x$ is
$$
u(x)= \mu(x)-\sum _{i=1}^{m} q_{i}(x)p_{i}[q(x),q(x)],
$$
whereas the utility (profit) value estimate of
$y$ is
$$
\varphi(x,y)= \mu(y)-\sum _{i=1}^{m} q_{i}(y)p_{i}[q(x),q(y)].
$$
Also, we set
$$
D(x)=\left\{ y \in X \ | \ q(y) \in U(q(x)) \right\}.
$$
Besides, we suppose that changing the activity profile
may invoke the necessity to change the treatment technology.
In particular, this may require new facilities, which were not used before.
These treatment change expenses for the transition $(q(x) \to q(y))$ can be determined by the bi-function
$c[q(x),q(y)]$. Hence, we can define the estimate of
the pure expenses for the move $(x \to y)$ as follows
$$
f(x,y)=\varphi(x,x)+c[q(x),q(y)]-\varphi(x,y),
$$
which coincides with that in Section \ref{sc:2}.
Given the initial activity profile $x^{0} \in X$, Problem  {\bf (P2)} will consist in finding
 a feasible trajectory approximating a solution of Problem  {\bf (P1)}.
In such a way, one finds a  relatively optimal technology activity profile.}
\end{example}


\begin{example} \label{ex:5.2} {\bf (Resource allocation in
telecommunication networks).} {\em We first describe an optimal flow
distribution problem in telecommunication data transmission networks.
The network contains $n$
transmission links (arcs) and accomplishes some submitted data
transmission requirements from $n$ selected pairs of
origin-destination vertices within a fixed time period.
Denote by $z_{i}$ and $d_{i}$ the
current and maximal value of data transmission for pair demand $i$,
respectively, and by $x_{j}$ the capacity of link $j$. Each pair
demand is associated with a
unique data transmission path, hence each link $j$ is associated uniquely with
the set $N(j)$ of pairs of origin-destination vertices, whose
transmission paths contain this link. For each pair demand $i$ we
denote by $\mu_{i}(z_{i})$ the network profit value at the data
transmission volume $z_{i}$. Then we can write the network profit
maximization problem as follows:
$$
\max  \to  \mu(z)=\sum \limits_{i=1}^{m} \mu_{i}(z_{i})
$$
subject to
\begin{eqnarray*}
  && \sum \limits_{i \in N(j)} z_{i} \leq  x_{j}, \ j=1, \dots, n;  \\
 && 0 \leq z_{i} \leq d_{i}, \ i =1, \dots, m.
\end{eqnarray*}

Denote by $u(x)$ the optimal value of this problem depending on
the right-hand sides $x$ of the constraints as parameters.
Let  $X$ denote the set of all the feasible capacity profiles, for instance, we can take
$$
X = \left\{x \in \mathbb{R}^{n} \ \vrule \  0 \leq x_{j} \leq
\alpha_{j}, \ j =1, \dots, n, \  \sum \limits_{j=1}^{n} \beta_{j} x_{j} \leq
C \right\}.
$$
That is,  $X$ stands for the set of feasible states. Each capacity profile
$x$ reflects the fixed allocation of network resources, hence, the transition
$(x \to y)$ requires certain expenses $c(x,y)$.
Suppose that one can calculate the values $c(x,y)$ and $u(y)$
only if $y$ belongs to some neighborhood $D_{1}(x)$ of $x$ and
that the direct transition $(x \to y)$ is possible  within the fixed time period
only if $y$ belongs to some neighborhood $D_{2}(x)$ of $x$.
In fact, some deviations from the current capacity profile may require new facilities, which were not
used before and essential changes in network organization.
Then we can set $D(x)=D_{1}(x)\bigcap D_{2}(x)$.
Given a current state $x^{0} \in X$,  Problem  {\bf (P2)} will determine
 a feasible trajectory of allocations tending to a relatively optimal solution.}
\end{example}


\section{Conclusions}

We presented  a rather general class of relative
optimization problems in metric spaces. The stationary
problem is formulated  as a quasi-equilibrium problem since
the goal function and feasible set are dependent of states.
The dynamic problem consists in finding a trajectory
attributed to an initial state  such that its points tend to
a solution of the stationary problem.
We proposed simple descent solution methods for creating
suitable trajectories to a relatively optimal state
under different conditions, which also gave
existence results for these problems. The approach was
illustrated by applied models.


\end{document}